\documentclass[12pt]{article}
\usepackage{tikz}
\usepackage{microtype}
\usepackage{geometry}                
\usepackage{amssymb,amsmath,amsthm}
\newtheorem{theorem}{Theorem}
\newtheorem*{theorems}{Theorems}
\theoremstyle{remark}
\newtheorem*{remark}{Remark}
\usepackage{color}
\usepackage{hyperref}
\usepackage{wrapfig}
\usepackage{graphicx}
\usepackage{url}

\usepackage{ifxetex}
\ifxetex
\usepackage{fontspec,xltxtra,xunicode}
\usepackage{unicode-math}
\defaultfontfeatures{Mapping=tex-text}
\newfontfamily\cyrillicfont{STIX Two Text}
\setromanfont[Mapping=tex-text]{STIX Two Text}
\usepackage{polyglossia}
\setotherlanguage{russian}
\setotherlanguage{german}
\setdefaultlanguage{english}
\else
\usepackage{mathptmx}
\fi


\let\le=\leqslant

\begin{document}
\title{Priority arguments and separation problems\footnote{Doklady AN SSSR (Soviet Math. Dokl.), volume 248, issue 6 (1979), pages 1309--1313. Translated by the author (2018)}}
\author{Alexander Shen\thanks{Recommended by A.N.~Kolmogorov, March 13, 1979}}
\date{}

\maketitle

Different constructions in the recursion theory use the so-called priority arguments (see, e.g.,~\cite[Section 13, The Priority Method]{shoenfield}). A general scheme was suggested by A.~Lachlan in~\cite{lachlan}. Based on this paper, we define the notion of a priority-closed class of requirements. Then, for a specific priority construction, we need to check only that all requirements we want to satisfy belong to some priority-closed class (defined in game terms). This game version of Lachlan's approach is used to present some results about recursively inseparable sets obtained by the author.

\fbox{\textbf{1}} A \emph{framework} is a quadruple $S=\langle M, U, R, T\rangle$, where $M$ and $U$ are some sets, and $R\subset M\times M$ and $T\subset M\times U$ are some relations. We assume that $T$ is not empty. The elements of $M$ are called \emph{items}, the elements of $U$ are called \emph{requests}. We read $R(m_1,m_2)$ as ``the element $m_1$ may follow the element $m_2$'', and $T(m,u)$ as ``the element $m$ is compatible with the request $u$''.

Here are two examples of frameworks. 
\begin{enumerate}
\item Items are finite subsets of $\mathbb{N}$ (the set of natural numbers), requests are pairs of disjoint finite subsets of $\mathbb{N}$.  We let
\[
R(A,B):= (B\subset A); \quad
T(A, \langle A^+, A^-\rangle):= (A^+\subset A\subset \mathbb{N}\setminus A^-).
\]
This framework will be denoted by $S_{+-}$.

\item Let $a(n)$ be a computable enumeration without repetitions of an enumerable [=c.e.] undecidable [=non-computable] set $A$. Let items be finite bit sequences; $R(m,n)$ means that $n$ is a proper prefix of $m$. For each item $m=m(0)\ldots m(k)$ consider the sets $A^0(m)$ and $A^1(m)$ defined as
\[
A^i(m)=\{ a(s) \mid m(s)=i\}.
\]
The requests are pairs of disjoint finite sets $\langle C_0, C_1\rangle$;
\[
T(m,\langle C_0,C_1\rangle):= (A^0(m)\cap C^0=A^1(m)\cap C^1=\varnothing).
\] 
This framework is denoted by $S^A_{01}$ in the sequel.
\end{enumerate}

\fbox{\textbf{2}} Let $S$ be some framework. We say that a request $u$ \emph{dominates} a request $v$ in $S$ if every item that is compatible with $u$ is also compatible with $v$. A sequence $m_0m_1\ldots$ of items (of $S$) is \emph{valid} if $m_{i+1}$ may follow $m_i$ for every $i$. A set $\alpha$ of valid sequences is called a \emph{condition} if it is closed under adding a prefix that keeps the sequence valid; if a sequence $m$ belongs to $\alpha$, we say that ``$m$ satisfies $\alpha$''. A \emph{requirement} is a pair of the form $\langle S,\alpha\rangle$ where $S$ is a framework and $\alpha$ is a condition (for $S$). A requirement $\langle S,\alpha\rangle$ is \emph{weaker} than a requirement $\langle S,\beta\rangle$ if $\beta\subset \alpha$. The \emph{conjunction} of two requirements $\langle S,\alpha\rangle$ and $\langle S, \beta\rangle$ is the requirement $\langle S, \alpha\cap\beta\rangle$.

Let $S_1$ and $S_2$ be two frameworks. We define the framework $S_1\times S_2$ as follows. It items and requests are pairs of items (resp. requests) for $S_1$ and $S_2$, the relations $R$ and $T$ and defined component-wise. Let $\alpha$ be a condition for $S_1$; we define the condition $\alpha\times S_2$ for $S_1\times S_2$ that is satisfied by sequences of pairs whose $S_1$-projection (the sequence of first components) satisfies $\alpha$. In a similar way we define the product of $S_1$ and some condition for $S_2$.

A class of requirements is \emph{priority-closed} if the following is true:
\begin{enumerate}
\item a requirement that is weaker than some requirement in the class, is also in the class;
\item the conjunction of any two requirements in the class is also in the class;
\item if a requirement $\langle S,\alpha\rangle$ is in the class and $S'$ is a framework, then the requirements $\langle S,\alpha\rangle\times S'$ and $S'\times \langle S,\alpha\rangle$ are in the class;
\item if a requirement $\langle S,\alpha\rangle$ is in the class, then $\alpha\ne\varnothing$.
\end{enumerate}

\fbox{\textbf{3}} For every requirement $\langle S,\alpha\rangle$ we consider a two-player game [with full information]. The players are Alice and Bob. Alice starts the game by choosing a request $u_0$ and an item $m_0$ that is compatible with $u_0$ (in the framework $S$). Then Bob chooses some request $u_1$ that dominates $u_0$. At the next move Alice chooses some item $m_1$ that may follow $m_0$ and is compatible with $u_1$. Then Bob chooses some request $u_2$ that dominates $u_0$ (but may not dominate $u_1$), Alice chooses some item $m_2$ that may follow $m_1$ and is compatible with $u_2$, and so on. The game is infinite. Bob wins in the game if
\begin{enumerate}
\item the game is infinite (Alice always has a move that does not violate the rules of the game);
\item the valid sequence $\alpha$ of items that appear during the game satisfies the condition $\alpha$;
\item all Bob's requests are the same starting from some moment in the game.
\end{enumerate}

A requirement $\langle S,\alpha\rangle$ is called a \emph{priority requirement} if Bob has a winning strategy in the corresponding game. Let $K_\textrm{p}$ be the class of all priority requirements. A requirement $\langle S,\alpha\rangle$ is called a \emph{countable-priority} requirement if it is the intersection of a countable family of priority requirements $\langle S,\alpha_0\rangle, \langle S,\alpha_1\rangle,\ldots$ that share the same framework $S$. Let $K_{\omega\mathrm{p}}$ is the class of all countable-priority requirements.

A framework $S=\langle M,U,R,T\rangle$ is \emph{constructive} if $M$ and $U$ are spaces of constructive finite objects (see~\cite{shoenfield}) and the relations $R$ and $T$ are decidable [=computable].  A \emph{computable priority requirement} has a constructive framework, and the corresponding game has a computable winning strategy for Bob. A \emph{weak computable priority requirement} has a constructive framework and a computable strategy for Bob that lets Bob win against computable strategies of Alice (=in all games where the sequence of items is computable and Bob follows the strategy). The class of all computable priority requirements is denoted by $K_{\mathrm{cp}}$, the class of all weakly computable priority requirement is denoted by $K_{\mathrm{wcp}}$. A requirement $\langle S,\alpha\rangle$ is \emph{countably computable priority requirement} if the framework $S$ is constructive and there exist computable priority requirements $\langle S,\alpha_0\rangle, \langle S, \alpha_1\rangle,\ldots$ such that $\alpha=\bigcap \alpha_i$ and the winning strategy for Bob in the $i$th game is computable uniformly in~$i$. In a similar way the \emph{countably weakly computable priority} requirements are defined. These two classes are denoted by $K_{\omega\mathrm{cp}}$ and $K_{\omega\mathrm{wcp}}$. The following inclusions immediately follow from the definitions:
\[
K_{\mathrm{p}}\subset K_{\omega\mathrm{p}},
K_{\mathrm{wp}}\subset K_{\omega\mathrm{wp}},
K_{\mathrm{cwp}}\subset K_{\omega\mathrm{cwp}},
K_{\mathrm{wp}}\subset K_{\mathrm{p}},
K_{\omega\mathrm{wp}}\subset K_{\omega\mathrm{p}},
K_{\mathrm{cp}}\subset K_{\mathrm{wcp}},
K_{\omega\mathrm{wp}}\subset K_{\omega\mathrm{wcp}}.
\]

\begin{theorem}\label{th:main}
The classes $K_{\mathrm{p}}$, $K_{\omega\mathrm{p}}$, $K_{\omega\mathrm{wcp}}$ are priority closed. If $\langle S,\alpha\rangle\in K_{\omega\mathrm{wcp}}$, then there exists a computable sequence that satisfies $\alpha$.
\end{theorem}

The proofs are omitted due to the lack of space.

\begin{remark}
The statement remains true if we replace $K_{\omega\mathrm{wcp}}$ by $K_{\mathrm{cp}}$, $K_{\omega\mathrm{cp}}$, or $K_{\mathrm{wcp}}$.
\end{remark}

\fbox{\textbf{4}} The proofs of several known existence results in the computability theory can be presented in the following form. We define some requirements and show that they belong to one of these classes. Then we consider some combination of these requirements that also belongs to the same class since the class is priority-closed. Then we use Theorem~\ref{th:main} and get a computable sequence that satisfies all the requirements and provides an object we are looking for. In the simple cases the class $K_{\omega\mathrm{cp}}$ is enough, but sometimes the larger class $K_{\omega\mathrm{wcp}}$ is needed. 

For every valid sequence $A_0\subset A_1\subset\ldots$ of items in the framework $S_{+-}$ (see section \textbf{1}) we consider the set $A_{\infty}=\bigcup A_i$. Let $\alpha$ be one of the conditions (1)~``the complement of $A_{\infty}$ is hyperimmune'', (2)~``the complement of $A_{\infty}$ is not hyperhyperimmune''. Then $\langle S_{+-},\alpha\rangle$ belongs to $K_{\omega\mathrm{cp}}$. Let $C$ be an undecidable set such that $C\le_T \mathbf{0}'$ and let $\alpha$ be the condition ($3{_C}$) ``$C$ is not reducible to $A_{\infty}$'' (here and below we use Turing reducibility), then $\langle S_{+-},\alpha\rangle$ belongs to $K_{\omega\mathrm{wcp}}$. Theorem~\ref{th:main} now guarantees that there exist an enumerable [=c.e.] set that satisfies any combination of the conditions of type (1), (2), ($3_C$).

In the framework $S_{+-}\times S_{+-}$ each valid sequence of items $\langle A_0,B_0\rangle, \langle A_1,B_1\rangle,\ldots$ determines two sets $A_{\infty}$ and $B_{\infty}$, the unions of $A_i$ and $B_i$. Let $\alpha$ be the condition ``$A_{\infty}$ is Turing incomparable with $B_{\infty}$''. One can show that $\langle S_{+-}\times S_{+-},\alpha\rangle$ belongs to $K_{\omega\mathrm{cp}}$. This implies Friedberg--Muchnik result~\cite{rogers} saying that there exist two incomparable enumerable sets. One can also ensure that these sets have additional properties (are hypersimple, not hyperhypersimple, etc.)

In the framework $S^A_{01}$ each valid sequence of items defines an infinite bit sequence $m(0)m(1)\ldots$ (that extends all items) and a splitting of the set $A$ into two sets $A^0$ and $A^1$ defined as $A^i=\{a(s)\mid m(s)=i\}$. Let $\alpha$ be one of the conditions: ``$A^0$ and $A^1$ are inseparable'', ``$C$ is not reducible to $A^0$'', ``$C$ is not reducible to $A^1$'' (for any undecidable set $C\le_T\mathbf{0}'$), ``$A^0$ and $A^1$ are incomparable''. Then one can show that $\langle S^{A}_{01},\alpha\rangle$ belongs to $K_{\omega\mathrm{wcp}}$. After that Theorem~\ref{th:main} can be applied to derive several corollaries, including the following known results.
\begin{theorems}
$1^*$.~\textup{\cite{friedberg}} Every enumerable undecidable set can be represented as the union of two disjoint enumerable undecidable sets.

$2^*$.~\textup(G.E.~Sacks, see~\textup{\cite{shoenfield}}\textup) Every enumerable undecidable set can be represented as the union of two disjoint incomparable enumerable sets.
\end{theorems}

\fbox{\textbf{5}} Consider the connections between priority arguments and diagonal arguments. A pair $\langle M, R\rangle$ is a \emph{$D$-framework} if $M$ is a set and $R\subset M\times M$ is a binary relation on $M$. The notions of a valid sequence and a condition are defined as before in section~\textbf{1}. A pair of the type $\langle S,\alpha\rangle$ where $S$ is a $D$-framework and $\alpha$ is a condition in $S$, is called \emph{$D$-requirement}. A $D$-requirement $\langle S,\alpha\rangle$ is called a \emph{diagonal} $D$-requirement if for every item $m_0$ of $S$ there exists a finite valid sequence that starts with $m_0$ such that every its infinite extension satisfies $\alpha$. 

\begin{theorem}\label{th:diagonal}
For every framework $S=\langle M,U,R,T\rangle$ there exists a $D$-framework $S_1=\langle M_1,R_1\rangle$ and a mapping $\varphi$ of the set of all valid sequences of $S_1$ into the set of all valid sequences of $S$ such that for every condition $\alpha$ in $S$ the following holds: if $\langle S,\alpha\rangle$ is a priority requirement, then $\langle S_1,\varphi^{-1}(\alpha)\rangle$ is a diagonal $D$-requirement.
\end{theorem}

In this theorem one can also exchange  ``framework'' and ``$D$-framework'', ``priority'' and ``diagonal'', ``requirement'' and ``$D$-requirement'' (simultaneously). Therefore, the use of classes $K_{\mathrm{p}}$ and $K_{\omega\mathrm{p}}$ is equivalent to diagonal arguments.

\fbox{\textbf{6}} Consider a framework $S^{++}_{-}$ that will be used to construct pairs of disjoint enumerable sets with required properties. Let items be the pairs of disjoint finite sets; the item $\langle A_1,B_1\rangle$ may follow $\langle A,B\rangle$ if $A\subset A_1$ and $B\subset B_1$. Let the requests be triples $\langle A^+, B^+, K\rangle$ of pairwise disjoint finite sets. An item $\langle A,B\rangle$ is compatible with a requirement $\langle A^+, B^+, K\rangle$ if $A^+\subset A$, $B^+\subset B$ and $(A\cup B)\cap K=\varnothing$.

For every valid sequence $\langle A_0,B_0\rangle, \langle A_1,B_1\rangle,\ldots$ of items in this framework we consider two sets $A_{\infty}=\bigcup A_i$ and $B_{\infty}=\bigcup B_i$. Let $\alpha$ be one of the following conditions: ``$A_{\infty}$ and $B_{\infty}$ and strongly inseparable''~\cite[exercise 8.39, p. 125]{rogers}, ``$A_{\infty}$ and $B_{\infty}$ are incomparable''. Then one can prove that the requirement $\langle S^{++}_{-},\alpha\rangle$ belongs to $K_{\omega\mathrm{cp}}$. 

We say that a set $D$ is a \emph{separator} for a pair $\langle X,Y\rangle$ of disjoint sets if $X\subset D$ and $Y\cap D=\varnothing$. One can prove that if $\alpha$ is one of the following conditions: ``no separator of $\langle C_1,C_2\rangle$ is reducible to any separator of $A_{\infty}, B_{\infty}$'' (for a given pair $\langle C_1,C_2\rangle$ of enumerable inseparable sets); ``$C$ is not reducible to any separator of $\langle A_{\infty},B_{\infty}\rangle$'' (for any undecidable $C\le_T\mathbf{0}'$), then the requrement $\langle S^{++}_{-},\alpha\rangle$ belongs to $K_{\omega\mathrm{wcp}}$. 

In the framework $S^{++}_{-}\times S^{++}_{-}$ every valid sequence of items determines two pairs $\langle A_{\infty},B_{\infty}\rangle$ and $\langle C_{\infty},D_{\infty}\rangle$. One can prove that for the condition $\alpha$=``any separator of $\langle A_{\infty},B_{\infty}\rangle$ is incomparable with any separator of $\langle C_{\infty}, D_{\infty}\rangle$'' the requirement $\langle S^{++}_{-}\times S^{++}_{-},\alpha\rangle$ belongs to $K_{\omega\mathrm{cp}}$. 

Theorem~\ref{th:main} then can be used to get the following results as corollaries:

\begin{theorem}\label{th:separators}
There exist two pairs $\langle A_1,A_2\rangle$ and $\langle B_1,B_2\rangle$ of enumerable inseparable sets such that any separator of $\langle A_1,A_2\rangle$ is incomparable with any separator of $\langle B_1,B_2\rangle$. Both pairs could be made strongly inseparable. One may also require that $A_1$ and $A_2$ are incomparable and also $B_1$ and $B_2$ are incomparable.
\end{theorem}

\begin{theorem}
For every pair $\langle C_1,C_2\rangle$ of enumerable inseparable sets there exists a pair $\langle A_1,A_2\rangle$ of enumerable inseparable sets such that no separator of $\langle C_1,C_2\rangle$ is reducible to any separator of $\langle A_1,A_2\rangle$. One may also require that the pair $\langle A_1,A_2\rangle$ is strongly inseparable.
\end{theorem}

\begin{theorem}
For every undecidable set $C\le_T\mathbf{0}'$ there exists a pair $\langle A_1,A_2\rangle$ of enumerable inseparable sets such that $C$ is not reducible to any separator of $\langle A_1,A_2\rangle$. 
\end{theorem}

\fbox{\textbf{7}} Consider one more framework that can be used to study separation problems. Let $\langle A_1,A_2\rangle$ be a pair of enumerable inseparable sets, and let $a(n)$ be a computable enumeration of $A_1$ without repetitions. Consider finite sets as items and pairs $\langle K^+,A^-\rangle$ such that $a(K^+)\cap A^-=\varnothing$ as requests. We say that item $K_2$ may follow $K_1$ if $K_1\subset K_2$. We say that item $K$ is compatible with request $\langle K^+,A^-\rangle$ if $K^+\subset K$ and $a(K)\cap A^-=\varnothing$.  We denote this framework by $S_{\oplus-}$. For every valid sequence of items $K_0\subset K_1\subset\ldots$ consider the set $A_1'=a(\bigcup K_i)$. One can prove that if $\alpha$ is one of the conditions ``$A_1'$ and $A_2$ are inseparable'', ``$C$ is not reducible to $A_1'$'' (for any undecidable $C\le_T\mathbf{0}'$), ``no separator for $\langle C_1,C_2\rangle$ is reducible to $A_1'$'' (for any pair $\langle C_1,C_2\rangle$ of enumerable inseparable sets), then $\langle S_{\oplus -},\alpha\rangle$ belongs to $K_{\omega\mathrm{wcp}}$. This implies, in particular, the following result.

\begin{theorem}
For every unsolvable separation problem in the sense of Yu.T.~Medvedev\textup{\cite{medvedev}} there exists a strictly smaller with respect to weak \textup[=Muchnik\textup] reducibility \textup[unsolvable\textup] separation problem.
\end{theorem}

The definition of weak reducibility can be found in~\cite{muchnik}.

\bigskip
\hbox to \textwidth{%
\hbox{\small Moscow Lomonosov State University}
\hfill
\hbox{\small Received April 20, 1979}
}
\par

\end{document}